\documentclass[reqno]{amsart}
\usepackage{amssymb}

\def\Var{\mathrm{Var}}
\def\idd{\mathop{\fam 0 id}\nolimits}
\def\charact{\mathop{\fam 0 char}\nolimits}

\def\Span{\mathop{\fam 0 Span}\nolimits}

\def\Cend{\mathop{\fam 0 Cend}\nolimits}
\def\Chom{\mathop{\fam 0 Chom}\nolimits}
\def\End{\mathop{\fam 0 End}\nolimits}
\def\gc{\mathop{\fam 0 gc}\nolimits}
\def\Curr{\mathop{\fam 0 Cur}\nolimits}
\def\alg{\mathrm{alg}}

\def\oo#1{\mathrel{{}_{(#1)}}}

\begin{document}

\title{Conformal representations of Leibniz algebras}
\thanks{Partially supported by RFBR.
The author gratefully acknowledges the support from Deligne
2004 Balzan
prize in mathematics.}
\author{Pavel Kolesnikov}
\address{Sobolev Institute of Mathematics, Novosibirsk, Russia}
\email{pavelsk@math.nsc.ru}
\begin{abstract}
In this note we present a more detailed
and explicit exposition of the
definition of a conformal representation of a Leibniz algebra.
Recall  (arXiv:math/0611501v3) that Leibniz algebras are exactly Lie
dialgebras.
The idea is based on the general fact that every dialgebra that
belongs to a variety $\Var $ can be embedded into
a conformal algebra of the same variety.
In particular, we prove that an arbitrary (finite dimensional)
Leibniz algebra has a (finite) faithful conformal representation.
As a corollary, we deduce the analogue of the PBW-theorem
for Leibniz algebras.
\end{abstract}

\maketitle

{\bf Introduction.}
In [1] we studied the relation between conformal algebras
[2] (these are algebraic systems with a countable number
of multiplications
coming from mathematical physics~[3])
and dialgebras [4,\,5]
(these are systems with two operations of multiplication
that were introduced in relation with Leibniz algebras).
We have proposed a general notion of what is a variety  of dialgebras
generalizing the classes of associative~[4] and alternative~[6]
dialgebras.
In particular, Lie dialgebras are exactly Leibniz algebras.
The main result of~[1] states that an arbitrary dialgebra of a given
variety can be embedded (in a canonical way) into a conformal
algebra of
the same variety.

In this paper, we focus on the embedding of Lie dialgebras (Leibniz
algebras) into conformal algebras.
For this variety, the construction of a conformal envelope from~[1]
can be simplified by using the notion of a conformal representation.
We will show that every (finite dimensional) Leibniz algebra can
be embedded into current conformal algebra over the algebra of linear
transformations of a (finite dimensional) linear space.
As a corollary, we obtain a new proof of the theorem on injective
embedding of a Leibniz algebra into an associative dialgebra and,
more explicitly, prove an analogue of the Poincare---Birkhoff---Witt
(PBW) theorem for Leibniz algebras~[5].

\medskip
{\bf 1. Varieties of dialgebras.}
A dialgebra is a linear space
$A$ endowed with two bilinear operations
$\dashv$, $\vdash $. In [4] and [6], the notions of
associative and alternative dialgebras were introduced.
These definitions were motivated by the relations between dialgebras
and Leibniz algebras. In particular, the  associativity identities
for dialgebras are chosen in such a way that the new operation
$$
  [ab] = a\vdash b - b\dashv a
$$
makes a dialgebra $A$ to be a left Leibniz algebra, i.e.,
the identity
$$
 [x_1[x_2x_3]] = [[x_1x_2]x_3] + [x_2[x_1x_3]]
\eqno{(1)}
$$
holds.
The Leibniz algebra $A^{(-)}$ obtained is an analogue of
the adjoint Lie algebra of an associative algebra.

In [1], we developed a general approach to the definition
of what is a variety of dialgebras via the notion of an operad.
In terms of identities, this definition can be stated as follows.

\smallskip
\noindent{\sc Definition 1.}
Let $\Var $
be a variety of algebras defined by a family
$\Sigma $
of polylinear identities.
A dialgebra $A$ is said to be a
$\Var$-dialgebra,
if it satisfies the identities
$$
\begin{gathered}
(x_1\dashv x_2)\vdash x_3 - (x_1\vdash x_2)\vdash x_3,
\quad
x_1 \dashv (x_2\vdash x_3) - x_1 \dashv (x_2\dashv x_3) , \\
t_i(x_1, \dots, x_n), \quad t\in \Sigma,\ \deg t=n, \ i=1,\dots, n,
\end{gathered}
 \eqno{(2)}
$$
where  $t_i$ is obtained from $t$
by replacing of each monomial
$(x_{j_1}\dots x_{j_n})$
with
$(x_{j_1}\vdash \dots \vdash x_i \dashv \dots \dashv x_{j_n})$
preserving the bracketing.
\smallskip

If $\Sigma = \{(x_1x_2)x_3 - x_1(x_2x_3) \}$
(associativity) then the system of identities~(2)
coincides with the one introduced in~[4]:
$$
\begin{gathered}
 t_1 = (x_1\dashv x_2)\dashv x_3 - x_1\dashv (x_2\dashv x_3),
\quad
 t_2 = (x_1\vdash x_2)\dashv x_3 - x_1\vdash (x_2\dashv x_3), \\
 t_3 = (x_1\vdash x_2)\vdash x_3 - x_1\vdash (x_2\vdash x_3).
\end{gathered}
$$
If
$\Sigma = \{x_1x_2+x_2x_1, \, x_1(x_2x_3) - x_2(x_1x_3)-(x_1x_2)x_3 \}$
(Lie algebra identities)
then, with respect to the operation
$[xy]=x\vdash y$,
the identities~(2) are equivalent to the left Leibniz identity~(1).
Conversely, an arbitrary (left) Leibniz algebra with respect to
$x\vdash y = [xy]$, $x\dashv y = -[yx]$
is a Lie dialgebra in the sense by Definition~1.

\medskip
{\bf 2. Conformal algebras.}
Let $H$ be a cocommutative Hopf algebra
with coproduct $\Delta$, counit $\varepsilon $
and antipode~$S$.
We will use the Sweedler notation~[7] omitting the summation symbols,
e.g.,
$$
\begin{gathered}
\Delta(h)= h_{(1)}\otimes h_{(2)}, \quad (\Delta \otimes
\idd)\Delta(h) = h_{(1)}\otimes h_{(2)} \otimes h_{(3)}, \\
(S \otimes
\idd)\Delta(h) = h_{(-1)}\otimes h_{(2)},
\quad
(\idd \otimes S)\Delta(h) = h_{(1)}\otimes h_{(-2)}. \\
\end{gathered}
$$

\noindent
{\sc Definition 2 {\rm [2,\,8,\,9]}.}
Let $M$ and  $N$ be two left $H$-modules.
A conformal linear map from
$M$ to~$N$
is a linear map  $a$
from $M$ to $H^{\otimes 2}\otimes_H N$
such that
$a(hu) = (1\otimes h\otimes_H \idd)a(u)$,
$h\in H$, $u\in M$.
Hereinafter, $H^{\otimes 2}= H\otimes H$
is considered as the outer product of regular right $H$-modules.
\smallskip

Note that $H\otimes H$ is a free right $H$-module with a basis
$\{h_i\otimes 1\}$, where
$\{h_i\}$ is a linear basis of  $H$ over $\Bbbk $.
Indeed,
$f\otimes g = (fg_{(-1)}\otimes 1 )g_{(2)}
=  \sum\limits_{i} (h_i\otimes 1)g_i$,
and such a presentation is unique.

The space of all conformal linear maps from
$M$ to~$N$
is denoted by
$\Chom(M,N)$.
If
$N=M$ then we use
$\Cend M$ for  $\Chom (M,M)$.

Let us define the structure of a left $H$-module on
$\Chom (M,N)$
as follows:
$$
  (ha)(u) = (h\otimes 1\otimes_H \idd)a(u),
\quad h\in H,\ a\in \Chom(M,N),\ u\in M.
$$

\noindent
{\sc Definition 3 {\rm [9,\,10]}.}
A left $H$-module $C$ endowed with an $H$-linear map
$\mu: C\to \Cend C$
is called an $H$-conformal algebra (or pseudo-algebra over~$H$).
The operation
$(a*b) = \mu(a)b\in H^{\otimes 2}\otimes_H C$,
$a,b\in C$,
is called pseudo-product.
\smallskip

If $H=\Bbbk$ then the notion of an $H$-conformal algebra
coincides with the notion of an ordinary algebra over a field.
If $H=\Bbbk[T]$ (with respect to the canonical
Hopf algebra structure)
and $\charact \Bbbk=0$
then an $H$-conformal algebra is exactly what is called by conformal
algebra in~[2].
In this paper, we will preferably consider
$H$-conformal algebras over $H=\Bbbk[T]$
without a restriction on the characteristic of~$\Bbbk$.

The simplest example of an
$H$-conformal algebra is provided by the
construction of the current algebra
$\Curr A$
over an ordinary algebra~$A$.
Given an algebra $A$, one may construct
$\Curr A$
as the free $H$-module
$H\otimes A$
endowed with the $H$-linear map
$$
 \mu(f\otimes a): h\otimes b\mapsto f\otimes h\otimes_H (1\otimes ab), \quad a,b\in A,
\ f,h\in H.
$$

The pseudo-product
$*: C\otimes C \to  (H\otimes H)\otimes_H C$
on an $H$-conformal algebra $C$
is an $H^{\otimes 2}$-linear map.
A categorical interpretation of the notion of an algebra
over an operad implies the following generalization
(see~[10]):
the pseudo-product can be considered as a map
from
$(H^{\otimes n}\otimes_H C)\otimes (H^{\otimes m}\otimes_H C)$
to
$H^{\otimes (n+m)}\otimes_H C$
given by
$$
(F_1\otimes _H a_1)*(F_2\otimes _H a_2)
= (F_1\otimes F_2\otimes _H 1)(\Delta^{n-1}\otimes
\Delta^{m-1}\otimes_H \idd)
(a_1*a_2),
\eqno{(3)}
$$
where
$F_1\in H^{\otimes n}$,
$F_2\in H^{\otimes m}$,
$a_1,a_2\in C$,
$\Delta^0(h) = h$,
$\Delta^k(h)= h_{(1)}\otimes \ldots \otimes h_{(k+1)}$
for
$h\in H$,
$k\ge 1$.

Relation (3) allows to interpret an arbitrary polylinear term
$t$
(in the signature of an algebra over a field)
as a function on
$C^{\otimes n}$ ($n=\deg t$)
with values in
$H^{\otimes n}\otimes_H C$.

\smallskip
\noindent
{\sc Definition 4 {\rm [10]}.}
Let
$\Var $
be a variety of algebras defined by a family
$\Sigma $
of polylinear identities.
An
 $H$-conformal algebra $C$ over a cocommutative Hopf algebra $H$
belongs to a variety
$\Var$ (is a $\Var$-pseudo-algebra, for short)
if it satisfies the identities
$$
 \begin{gathered}
t^*(x_1, \dots, x_n), \quad t\in \Sigma,
 \end{gathered}
$$
where $t^*$
is obtained from
$t$
by replacing each monomial
$(x_{1\sigma}\dots x_{n\sigma})$,
$\sigma\in S_n$,
with
$(\sigma \otimes_H\idd)(x_{1\sigma }*\dots *x_{n\sigma })$
preserving the bracketing.
Here $ \sigma $ acts on $H^{\otimes n}$
by the rule
$(h_1\otimes\dots\otimes h_n)^\sigma  = (h_{1\sigma ^{-1}}\otimes
\dots \otimes h_{n\sigma ^{-1}})$.
\smallskip

\noindent
{\sc Remark.}
The class of $\Var $-pseudo-algebras over a Hopf algebra $H$
is not a variety in the ordinary sense: it is not closed under
infinite direct products.
\smallskip

On every $H$-conformal algebra $C$ one may define the
operations
$\dashv $, $\vdash $ by the
following rule:
if
$a*b = \sum\limits_i h_i\otimes 1\otimes_H c_i$,
$a,b,c_i\in C$,
$h_i\in H$,
then
$$
 a\vdash b =\sum\limits_{i} \varepsilon (h_i)c_i,
\quad
 a\dashv b  = \sum\limits_{i}h_ic_i.
$$
The dialgebra obtained is denoted by
$C^{(0)}$.

The main result of [1] can be stated in the following form.

\smallskip
\noindent
{\bf Theorem 1.}
{\sl
If $C$ is a $\Var $-pseudo-algebra
then
$C^{(0)}$
is a $\Var$-dialgebra.
Conversely, an arbitrary
$\Var$-dialgebra can be embedded into
$C^{(0)}$
for an appropriate
$\Var$-pseudo-algebra~$C$
over
$H=\Bbbk[T]$.}\quad $\square $

\medskip
{\bf 3. Faithful representations of Leibniz algebras.}
Let $H$ be the coordinate Hopf algebra
of an abelian linear algebraic group~$G$.
Let us fix a basis
$\{h_i\}$
of $H$ over~$\Bbbk$.
Consider three left
 $H$-modules $M$, $N$, and~$V$
together with an $H^{\otimes 2}$-linear map
$$
  *: M\otimes N \to H^{\otimes 2}\otimes _H V.
$$
Define the family of operations
$(\cdot\oo{z}\cdot): M\otimes N \to V$, $z\in G$,
as follows.
If
$u*v = \sum\limits_i
 h_i\otimes 1 \otimes _H w_i$
for $u \in M$,
$v\in N$,
then
$$
 (u\oo{z} v) = \sum\limits_{i} h_i(z^{-1})w_i \in V.
\eqno{(4)}
$$
The family of bilinear operations constructed
does not depend on the choice of
$\{h_i\}$, and it has the following properties:
\begin{itemize}
\item[(C1)]
for every $u$, $v$
the map
$x\mapsto (u\oo{x} v)$
is a regular
$V$-valued function on~$G$;
\item[(C2)]
$(hu\oo{z} v)=h(z^{-1})(u\oo{z}v)$;
\item[(C3)]
$(u\oo{z}hv)= h_{(1)}(z)h_{(2)}(u\oo{z} v)$.
\end{itemize}

The property
(C1) means that for every $u\in M$, $v\in N$
there exists a uniquely defined (if a basis $\{h_i\}$ is fixed)
finite set of $w_i\in V$
such that
$(u\oo{z} v) = \sum\limits_i h_i(z^{-1})w_i$
for all
$z\in G$.
In this settings, define
$\{u\oo{z}v\} = \sum\limits_i h_{i(1)}(z) h_{i(2)}w_i$.

The structure of an $H$-conformal algebra on a left
$H$-module $C$ is completely defined by the family of
operations~(4)
for
$M=N=V=C$.
The identities of conformal algebras described in Definition~4
can be expressed in terms of operations
$(\cdot\oo{z}\cdot)$
and $\{\cdot\oo{z}\cdot\}$.
For example,
an $H$-conformal algebra $C$ is associative
if and only if
$$
 a\oo{z} (b\oo{y} c) =
 (a\oo{z} b)\oo{yz} c,\quad a,b,c\in C,\ z,y\in G.
\eqno{(5)}
$$

\noindent
{\sc Remark.}
If
$G=\mathbb A_1\simeq (\Bbbk, +)$,
$\charact \Bbbk = 0$
then $(\cdot\oo{z}\cdot)$
coincides with what is called
$\lambda $-bracket in~[8].
\smallskip

If $M$
is a left $H$-module then the action of
$\Cend M$
on~$M$
is also completely defined by the family of
$G$-products
$(\cdot \oo{z} \cdot):  \Cend M\otimes M \to M$, $z\in G$,
built by the map
$a*u = a(u)$, $a\in \Cend M$, $u\in M$.
Define similar operations
$(\cdot\oo{z} \cdot)$
and
$\{\cdot \oo{z}\cdot \}$, $z\in G$,
on the space
$\Cend M$
by the following rule:
if
$a,b\in \Cend M$, $u\in M$,
$$
 b(u) = \sum\limits_i h_i\otimes 1\otimes_H v_i,
 \quad
 a(v_i) = \sum\limits_j h_j \otimes 1\otimes _H w_{ij},
$$
then
$$
 (a\oo{z} b)(u) = \sum\limits_{i,j}
    h_j(z^{-1})h_{i(1)}(z)h_{i(2)}\otimes
     1\otimes_H w_{ij}
\eqno{(6)}
$$
and
$$
 \{a\oo{z} b\}(u) = \sum\limits_{i,j} h_{i}(z)h_{j(1)}(z^{-1}) h_{j(2)}\otimes
     1\otimes_H w_{ij}.
\eqno{(7)}
$$

(It is easy to check that
 $(a\oo{z} b)$ and
$\{a\oo{z} b\} $
belongs to~$\Cend M$ for all $z\in G$.)
The operations~(6) on $\Cend M$
satisfy the relations
(C2),~(C3), and~ (5),
but (C1)
fails in general.

\smallskip
\noindent
{\bf Proposition 2.}
{\sl The linear space
 $\Cend M$
with respect to the operations
$a\vdash b = (a\oo{e} b)$,
$a\dashv b = \{a\oo{e} b\}$,
where $e$ is the unit of $G$,
is an associative dialgebra denoted by
$\Cend M^{(0)}$.}
\smallskip
\noindent
{\sc Proof.}
As we have already mentioned,
the operations (6) on $\Cend M$
satisfy~(5).
Moreover, the following identities hold:
\begin{gather}
\{\{a\oo{z} b\}\oo{y} c\} = \{a\oo{yz} \{b\oo{y} c\}\} = \{a\oo{yz} (b\oo{z^{-1}} c)\},
    \tag 8  \\
(\{a\oo{z} b\}\oo{y} c) = ((a\oo{yz} b)\oo{y} c) = (a\oo{yz} (b\oo{z^{-1}} c))
    \tag 9
\end{gather}
for all $y,z\in G$.

Indeed, the definitions (6) and (7) are chosen in such a way that
$$
 (a\oo{z} b)\oo{w} u = a\oo{z} (b\oo{wz^{-1} } u),
\quad
 \{a\oo{z}b\}\oo{w} u = a\oo{wz} (b\oo{z^{-1}} u)
$$
for $a,b\in \Cend M$, $u\in M$, $z,w\in G$.
These relations imply (8), (9), and also (5), e.g.,
$$
  \{\{ a\oo{z} b\}\oo{w} c\}\oo{y} u
= \{a\oo{z} b\}\oo{wy} (c\oo{w^{-1}} u)
= a\oo{wyz} (b\oo{z^{-1}} (c\oo{w^{-1}} u))
$$
and
$$
\{a\oo{wz} \{b\oo{w} c\}\}\oo{y} u
= a\oo{wyz} (\{b\oo{w} c\}\oo{(wz)^{-1}} u )
= a\oo{wyz} (b\oo{z^{-1}} (c\oo{w^{-1}} u)).
$$
Other equations in (8), (9) can be deduced analogously.

If $y=z=e$ then (5), (8), (9)
turn into the identities of an associative dialgebra.
\quad $\square $
\smallskip

The adjoint Lie dialgebra (Leibniz algebra)
with the operation
$[ab] = (a\oo{e} b) - \{b\oo{e} a\}$, $a,b\in \Cend M$,
is denoted by
$\gc M^{(0)}$.

\smallskip
\noindent
{\sc Definition 5.}
A conformal representation of an associative
(resp., Lie) dialgebra $A$
on a left $H$-module $M$ is a homomorphism of dialgebras
$\rho : A \to C^{(0)}$,
where $C=\Cend M$ (resp., $\gc M$).
If $M$ is a finitely generated $H$-module then
$\rho $ is said to be a finite representation.
\smallskip

In the case of ordinary algebras for
$G=\{e\}$ Definition~5 coincides with the usual definition
of what is a representation of an algebra on a linear space.
By the classical Poincare---Birkhoff---Witt and Ado---Iwasawa
theorems an arbitrary (finite dimensional) Lie algebra
has a (finite dimensional) faithful representation.
A Leibniz algebra which is not a Lie algebra obviously has
no faithful representation for $G=\{e\}$
(as well as for any finite group),
but it turns out that a faithful representation exists for
affine line $G=\mathbb A_1$.

\smallskip
\noindent
{\bf Theorem 3.}
{\sl
An arbitrary (finite dimensional) Leibniz algebra
has a (finite) faithful conformal representation for
 $G=\mathbb A_1$. }

\smallskip
\noindent
{\sc Proof.}
Let $L$ be a Leibniz algebra,
$L^{\alg} = L/\Span \{[ab]+[ba]\mid a,b\in L \}$
is the ordinary Lie algebra.
Denote the image of an element
$a\in L$ in $L^{\alg }$
by~$\bar a$.

Suppose we have fixed a non-trivial
 $L^{\alg}$-module
$V$, and consider the space
$M_0 = V\oplus (L\otimes V)$.

Recall that the coordinate Hopf algebra
of $G=\mathbb A_1$
is the polynomial algebra
 $H=\Bbbk[T]$,
where
$\Delta(h)=h(T\otimes 1+1\otimes T)$,
$\varepsilon (h)=h(0)$,
$S(h)=h(-T)$
for
$h\in H$.
We will construct a homomorphism
 $\rho : L\to \gc M^{(0)}$
for
$M=H\otimes M_0$.
To define a conformal linear map
$\rho(a)$, $a\in L$,
it is enough to fix the values of
$(\rho (a)\oo{z} u) \in M$, $z\in \Bbbk$, $u\in M_0$.
Set
$$
\begin{gathered}
 \rho (a)\oo{z} v = \bar av + z(a\otimes v), \quad v\in V,\\
 \rho (a)\oo{z} (b\otimes v) = b\otimes \bar av +
                  [ab]\otimes v,\quad b\in L, \ v\in V.
\end{gathered}
\eqno{(10)}
$$
Relations (10) define an injective $H$-linear map from
$L$
to $\Cend M$.
Let us show that $\rho $ is a conformal representation
of the Leibniz algebra~$L$.
It is enough to check that
$$
\rho([ab]) = (\rho(a)\oo{0}\rho(b)) - \{\rho(b)\oo{0}\rho(a) \}
$$
for all $a,b\in L$.
Indeed, using (5) and~(9) we obtain
\begin{multline}
((\rho(a)\oo{0}\rho(b)) - \{\rho(b)\oo{0}\rho(a) \})\oo{z} v
=
 \rho (a)\oo{0} (\rho (b)\oo{z} v) - \rho (b)\oo{z}(\rho (a)\oo{0} v ) \\
=
\rho (a)\oo{0} (\bar bv + z(b\otimes v))
  - \rho (b)\oo{z} (\bar av)                \\
=
\bar a\bar b v + z(b\otimes \bar a v) + z[ab]\otimes v
    - \bar b\bar a v - z b\otimes \bar av \\
=
[\bar a,\bar b] v + z[ab]\otimes v
=
\rho ([ab])\oo{z} v,
\nonumber
\end{multline}
\begin{multline}
((\rho(a)\oo{0}\rho(b)) - \{\rho(b)\oo{0}\rho(a) \})\oo{z} (c\otimes v) \\
=
\rho (a)\oo{0} (c\otimes \bar bv + [bc]\otimes v )
   -\rho (b)\oo{z} (c\otimes \bar av + [ac]\otimes v)                     \\
=
c\otimes [\bar a,\bar b]v + [[ab]c]\otimes v
=
\rho ([ab])\oo{z} (c\otimes v)
\nonumber
\end{multline}
for $c\in L$, $v\in V$.

If $\dim L <\infty $
then $\dim M_0$
can be chosen to be finite dimensional, and in this case
$\rho $
is a finite representation.
\quad $\square $
\smallskip

It is clear from (10)
that the images of
$a\in L$ under $\rho $ belong to
$\Curr \End M_0 \subset \Cend M$.
Indeed,
$\rho (a) = 1\otimes a_0 - T\otimes a_1$,
where
$$
\begin{gathered}
a_0(v) = \bar av, \quad
a_0(b\otimes v) =  b\otimes \bar a v + [ab]\otimes v,\\
a_1(v) = a\otimes v,\quad  a_1(b\otimes v)=0
\end{gathered}
$$
for all $b\in L$, $v\in V$.

\smallskip
\noindent
{\bf Corollary 4.}
{\sl
Every (finite dimensional) Leibniz algebra
can be embedded into current conformal algebra
over the ordinary algebra of linear transformations
of a (finite dimensional) linear space.}\quad $\square $

\smallskip
\noindent
{\bf Corollary 5 {\rm [5]}.}
{\sl For every Leibniz algebra $L$ there exists an
associative dialgebra $A$ such that
$L\subset A^{(-)}$.}\quad $\square $
\smallskip

In the last statement, just consider
$A=(\Curr\End M_0)^{(0)}$.

Let $L$ be a Leibniz algebra.
Denote by $U(L)$
its universal enveloping associative dialgebra built in
[4,\,5].
Recall that $U(L)$
is a quotient of the free associative dialgebra
$F(L)$
generated by the space~$L$
by the ideal generated in
$F(L)$
by all elements of the form
$ a\vdash  b - b\dashv a - {[ab]}$,
$a,b\in L$.
We will use
 $\tilde c$ to denote the image
of $c\in L\subset F(L)$ in $U(L)$.

It follows from the universal property of
$U(L)$
that an arbitrary conformal representation
$\rho : L \to \gc M^{(0)} = (\Cend M^{(0)})^{(-)}$
induces a representation
$U(\rho): U(L) \to \Cend M^{(0)}$.
of the associative dialgebra $U(L)$.
If
$M$
is the module from Theorem~3 built on the
$L^{\alg}$-module $V=U(L^{\alg})$
then let us denote the corresponding conformal representation
$L \to \gc M^{(0)}$ by $\rho_U$.

\smallskip
\noindent
{\bf Theorem 6.}
{\sl The representation
$U(\rho_U)$
is faithful. }

\smallskip
\noindent
{\sc Proof.}
We will prove the theorem without a reference to the results
of [4,\,5].

Let us fix a linear basis
$B$ of $L$ in the following form:
$B=\{a_i\mid i\in I\}\cup \{b_j\mid j\in J\}$,
where
$\{\bar a_i\mid i\in I\}$ is a basis of $L^{\alg}$,
$\{y_j\mid j\in J\}$ is a basis of the kernel
of the natural homomorphism
$L\to L^{\alg}$.
Suppose the set $I$  to be well-ordered.

\smallskip
\noindent
{\bf Lemma 7.}
{\sl The elements of the form
$$
\tilde c\dashv \tilde a_{i_1} \dashv  \dots \dashv \tilde a_{i_n} ,
\quad
c\in B,\ i_1\le \dots \le i_n,\ n\ge 0,
\eqno{(11)}
$$
span
$U(L)$.}

\smallskip
\noindent
{\sc Proof.}
Note that the bracketing in (11) is not essential since
the operation $\dashv $ is associative.

By making use of the identities of associative dialgebras
and the rewriting rule
$\tilde x\vdash \tilde y =\tilde y\dashv \tilde x +\widetilde{[xy]} $,
$x,y\in L$,
one may present an arbitrary element of
$U(L)$
as a linear combination of words
$$
 \tilde c_1\dashv \tilde c_2 \dashv \dots \dashv \tilde c_n ,
\quad c_i\in B.
$$
Since
$x\dashv \tilde b_j=0$
for any $x\in U(L)$,
we may assume
$c_2,\dots, c_n\in \{a_i\mid i\in I\}$.
It remains to note that if
 $i_1>i_2$
then
\begin{multline}
x\dashv (\tilde a_{i_1}\dashv \tilde a_{i_2})
 = x\dashv (\tilde a_{i_2}\dashv \tilde a_{i_1})
+ x\dashv (\tilde a_{i_1}\dashv \tilde a_{i_2} - \tilde a_{i_1}\vdash \tilde a_{i_2})
+ x\dashv \widetilde{[a_{i_1}a_{i_2}]} \\
= x\dashv (\tilde a_{i_2}\dashv \tilde a_{i_1}) +
 x\dashv \widetilde{[a_{i_1}a_{i_2}]}.
\nonumber
\end{multline}
Therefore, an arbitrary element of
$U(L)$
can be presented as a linear combination of~(11).
\quad $\square $
\smallskip

To complete the proof it is enough to show
that the images of the elements~(11) under $U(\rho_U)$
are linearly independent.
Suppose $u$ is a word of the form~(11).
Compute
$(U(\rho_U)(u)\oo{z} 1)\in M$
for
$1\in U(L^{\alg})$.
Since
$U(\rho_U)$ is a homomorphism of dialgebras,
$$
w=U(\rho_U)(u) =
  \{\tilde c\oo{0} \{\tilde a_{i_1} \oo{0}  \dots
\{\tilde a_{i_{n-1}} \oo{0} \tilde a_{i_n}\}\dots \}\}.
$$
Hence, relation (9) implies
$$
(w\oo{z} 1) =
(\tilde c\oo{z} (\tilde a_{i_1} \oo{0} ( \dots \oo{0} (\tilde a_{i_n}\oo{0} 1)
\dots ))),
$$
which is easy to compute by (10):
$$
 (w\oo{z} 1) = \tilde c\oo{z} (\bar a_{i_1}\dots\bar a_{i_n})
=\bar c\bar a_{i_1}\dots\bar a_{i_n} + z(c\otimes \bar a_{i_1}\dots\bar a_{i_n}).
$$
The coefficients at $z$ obtained are linearly independent in $M$
for different words of the form~(11).
Hence, the images of these words under $U(\rho_U)$
are also linearly independent.\quad $\square $
\smallskip

\noindent
{\bf Corollary 8 {\rm (PBW-Theorem for Leibniz algebras [5,\,11])}.}
{\sl The universal enveloping associative dialgebra
$U(L)$
for Leibniz algebra
$L$
is isomorphic to
$L\otimes U(L^{\alg})$
as a linear space.}\quad $\square $


\begin{thebibliography}{10}

\bibitem{1}
Kolesnikov P.~S.,
Varieties of dialgebras and conformal algebras,
arXiv:math/0611501v3

\bibitem{2}
Kac V.~G.,
Vertex Algebras for Beginners.  University Lecture Series, vol.~10,
AMS, Providence, RI, 1996.

\bibitem{3}
Belavin~A.~A., Polyakov~A.~M., Zamolodchikov~A.~B.,
Infinite conformal symmetry in two-di\-men\-sional quantum field theory,
Nuclear Phys. 241 (1984) 333--380.

\bibitem{4}
Loday~J.-L., Pirashvili~T.,
Universal enveloping algebras of Leibniz algebras and homology,
Math. Ann. 296 (1993) 139--158.

\bibitem{5}
Loday~J.-L.,
Dialgebras, in: Dialgebras and Related Operads.
Lecture Notes in Mathematics, vol.~1763.
Springer Verl., Berlin, 2001, pp.~7--66.

\bibitem{6}
Liu D.,
Steinberg---Leibniz algebras and superalgebras,
J. Algebra 283 (1) (2005) 199--221.

\bibitem{7}
Sweedler~M.,
Hopf Algebras,
New York, Benjamin, 1969.

\bibitem{8}
Kac V.~G.,
Formal distribution algebras and conformal algebras,
Proc. XII-th International Congress in Mathematical Physics (ICMP'97),
International Press, Cambridge, MA, 1999, pp.~80--97.

\bibitem{9}
Bakalov~B., D'Andrea~A.,  Kac~V.~G.,
Theory of finite pseudoalgebras,
Adv. Math. 162 (1) (2001) 1--140.

\bibitem{10}
Kolesnikov~P.~S.,
Identities of conformal algebras and pseudoalgebras,
Comm. Algebra 34 (6) (2006) 1965--1979.

\bibitem{11}
Aymon~M., Grivel~P.-P.,
Un th\'eor\`eme de Poincar\'e---Birkhoff---Witt pour les alg\`ebres de Leibniz,
Comm. Algebra 31 (2) (2003)  527--544.

\end{thebibliography}
\end{document}